\documentstyle{article}
\newcommand{\res}{\mathop{\hbox{\vrule height 7pt width .5pt depth 0pt
\vrule height .5pt width 6pt depth 0pt}}\nolimits}
\newcommand{\qed}{\thinspace\null\nobreak\hfill\hbox{\vbox{\kern-.2pt\hrule
height.2pt
depth.2pt\kern-.2pt\kern-.2pt \hbox to2.5mm{\kern-.2pt\vrule width.4pt
\kern-.2pt\raise2.5mm\vbox to.2pt{}\lower0pt\vtop to.2pt{}\hfil\kern-.2pt
\vrule width.4pt\kern-.2pt}\kern-.2pt\kern-.2pt\hrule height.2pt depth.2pt
\kern-.2pt}}\par\medbreak}

\newtheorem{definition}{Definition}[section]

\newtheorem{theorem}[definition]{Theorem}
\newtheorem{proposition}[definition]{Proposition}

\newtheorem{remark}[definition]{Remark}
\newtheorem{example}[definition]{Example}

\newcommand{\rr}{{\bf R}}
\newcommand{\rn}{\rr^n}

\newcommand{\NN}{{\bf N}}

\newcommand{\ZZ}{{\bf Z}}
\newcommand{\Zn}{{\bf Z}^n}

\newcommand{\Mn}{{\rm I\! M}^{n\times n}}        

\newcommand{\Aa}{{\cal A}}           
\newcommand{\Ao}{\Aa(\Omega)}

\newcommand{\dx}{\,dx}

\newcommand{\dy}{\,dy}

\newcommand{\ie}{; {\it i.e., }}
\newcommand{\e}{\varepsilon}

\newcommand{\lp}{{\rm L}^p}

\newcommand{\lu}{{\rm L}^1}

\newcommand{\wup}{{\rm W}^{1,p}}

\newcommand{\ld}{\rm LD}                 
\newcommand{\ldp}{{\rm LD}^p}                 
\newcommand{\wuu}{{\rm W}^{1,1}}

\newcommand{\wui}{{\rm W}^{1,\infty}}

\newcommand{\orm}{(\Omega;\rr^m)}
\newcommand{\orn}{(\Omega;\rr^n)}                  
\newcommand{\loc}{{\rm loc}}

\newcommand{\dist}{{\rm dist}\,}
\newcommand{\spt}{{\rm spt}\,}

\newcommand{\Hom}{{\rm hom}}     

\newcommand{\ave}{-\hskip -.38cm\int}

\newcommand{\Om}{\Omega}

\begin{document}
\vspace*{2cm}
\begin{center}
{\large \bf HOMOGENIZATION OF PERIODIC MULTI-DIMENSIONAL STRUCTURES: 
THE LINEARLY ELASTIC/PERFECTLY PLASTIC CASE}
\vspace{1cm}\\
{\sc Nadia Ansini and Fran\c{c}ois Bille Ebobisse}
\end{center}
\vspace{2.5cm}
{\bf Abstract.}
In this paper we study the asymptotic behaviour via $\Gamma$-convergence 
of some integral functionals $F_{\e}$ which model some multi-dimensional 
structures and depend explicitly on the linearized 
strain tensor.
The functionals $F_{\e}$ are defined in particular classes 
of functions with bounded deformation while the limit problem is set in the
usual framework of Sobolev spaces or $BD(\Om)$. We also construct an
example of such functionals showing
that under some special assumptions we can have non local effects.
\vspace{0.125cm}\\
\section{Introduction}
In recent years there has been an increasing interest in the description of 
media whose
microscopic behaviour takes into account lower dimensional or
multi-dimensional structures and can be modeled by suitable integral 
functionals with 
respect to periodic measures
(see  \cite {BCP}, \cite{Zh}, \cite{ABF}, \cite{BBS}, \cite{ABC}). 
Zhikov studied in \cite{Zh} the homogenization of functionals $F_{\e}$ 
defined as
$$
F_\e (u,\Om )=\int _\Om f \Bigl({x\over\e }, \nabla u\Bigr)d\mu _\e
$$
on  $C^\infty\orm$,  
where $\mu _\e $ is defined by $\mu _\e (B):=\e ^n\mu ({1\over\e }B)$
with $\mu $ a fixed $1$-periodic Radon measure and $f$ is a
Borel function $1$-periodic in the first variable
(see also  Braides and Chiad\`o Piat \cite {BCP} for the case $\mu=\chi_E$ with
$E$ periodic, and
Bouchitt\'e, Buttazzo and Seppecher \cite{BBS} for relaxation results
in the case of general $\mu$). On the other hand, 
following the approach of Ambrosio, Buttazzo and Fonseca \cite{ABF}, which is 
somehow complementary
to the ``smooth approach'' described above, Ansini, Braides and Chiad\`o Piat 
studied  in
\cite{ABC} the asymptotic behaviour of energy functionals concentrated
on periodic multi-dimen\-sional structures, of the form
$$
F_\e (u,\Om )=\int _\Om f\Bigl ({x\over\e },{{dDu}\over {d\mu
    _\e}}\Bigr )d\mu _\e \,.
$$
In this case the problem is set in
the framework of Sobolev spaces $\wup_{\mu_{\e}}\orm$
with respect to the measure $\mu _\e$ of \cite{ABF}. We recall that 
$\wup_{\mu_{\e}}\orm$ is the space of
functions $u\in \lp\orm$ whose distributional derivative is a
measure absolutely continuous with respect to $\mu _\e $ with
$p$-summable density $dDu/d\mu _\e$. 
 A homogenization theorem for
$F_\e$ has been proved under a standard growth condition of order
$p$ on $f$ and a notion of $p$-homogenizability introduced for the
measure $\mu$ (see \cite{ABC} Theorem 3.5).

In the context of linear elasticity or perfect
plasticity, in place of considering energies depending on the deformation 
gradient
$Du$, it is more appropriate to
consider energy functionals depending explicitly on the {\it linearized strain
tensor} $Eu$. Our goal in this paper is to  study the asymptotic behaviour 
of functionals of the type
$$
F_\e (u,\Om )=\int _\Om f\Bigl ({x\over\e},{{dEu}\over {d\mu
    _\e}}\Bigr )d\mu _\e 
$$
defined in a particular class of functions
with bounded deformation denoted by $\ldp_{\mu _\e }(\Om)$ (introduced in 
Section 3). More precisely, $\ldp_{\mu _\e }(\Om)$
is
the space of functions $u\in \lp\orn$, whose deformation
tensor $Eu$ is a measure absolutely continuous with respect to $\mu_\e $ 
with $p$-summable density $dEu/d\mu_\e $. 
Using both classical and fine properties of functions
with bounded deformation and 
the same assumptions as in \cite{ABC} with a modified 
definition of `$p$-homogenizable measure', we prove in the first part of
the paper, a homogenization
theorem (Theorem \ref{24.12}). Precisely, we show the
existence of the $\Gamma $-limit of the functionals $F_\e $ 
with respect to $\lp$-convergence in the
Sobolev space $\wup\orn$, and with respect to $\lu$-convergence in $BD(\Om)$ 
(the space of 
functions with bounded deformation in $\Om$, that is the space of functions 
$u\in \lu\orn$
whose deformation tensor $Eu$ is a Radon measure with finite total variation 
in $\Om$, see \cite{ACD}). We show that the
$\Gamma $-limit admits an integral representation

$$
F_{\Hom}(u,\Om)=\int_\Om f_{\Hom} (Eu)\dx
$$
in $\wup\orn$; moreover, if $f$ is convex then
$$
F_{\Hom}(u,\Om)=\int _\Om f_{\Hom} ({\cal E}u)\dx
+\int _\Om f_{\Hom} ^\infty \Bigl ({{dE^s u}\over {d\vert
  E^s u\vert }}\Bigr )d\vert E^s u\vert
$$ 
in $BD(\Om )$, where ${\cal E}u$ is the density of the 
absolutely continuous part and $E^s u$ is the singular part of $Eu$ with
respect to the Lebesgue measure;
$f_{\Hom}$ is described by an asymptotic formula and $f_{\Hom} ^\infty$ 
denotes the recession function of $f_{\Hom}$ (see (\ref{rec})).

In the second part of this paper we show that when the scaling argument 
leading to the functionals $F_{\e}$ does not apply, 
non local effects can arise. More precisely, we
consider functionals of the type
$$F^\gamma_\e(u, \Om )=\e^\gamma\int_\Om f\Bigl
({x\over\e},{dEu\over d\mu_\e}\Bigr)d\mu_\e\, ,
$$
which in the previous approach tend to the null functional when $\gamma >0$, 
and we construct an explicit example showing that, with a suitable
choice of $\gamma $, $\mu_\e$ and of the convergence with respect to which the
$\Gamma $-limit is computed, we have a limit functional of a non local nature.

\section{Notation and preliminaries}
In the sequel $\Mn$ stands for the
space of $n\times n$ matrices and $\Mn_{sym}$ for the space of $n\times n$
 symmetric matrices. The letter $c$ will stand for an arbitrary fixed  
strictly-positive constant independent of the parameters under consideration,
whose value may vary from line to line. The symbols $(\cdot,\cdot)$ and 
$|\cdot|$ stand for the Euclidean scalar product and the Euclidean norm. 
The  Hausdorff $k$-dimensional measure and the Lebesgue measure  
in $\rn$ are denoted by ${\cal H}^k$ and ${\cal L}^n$ respectively. 
We write $|E|$ for the Lebesgue measure ${\cal L}^n$ of $E$.
We recall that for any two vectors $a$ and $b$ in $\rn$, the symmetric
 product $a\odot b$ is the symmetric $n\times n$ matrix defined by 
$ a\odot b= {1 \over 2}(a \otimes b + b\otimes a)$, where $\otimes$ 
denotes the tensor product.
$\Omega$ is a bounded open subset of $\rn$; we denote
by $\Ao$ the family of all open subsets of $\Omega$.

Given a matrix-valued measure $\mu$ on $\Omega$, we adopt the
notation $|\mu|$ for its total variation (see Federer \cite{F}).
The measure $\mu \res F$ is defined by $(\mu \res F)(B)= \mu(B\cap F)$.
We write $\mu <\!<\lambda$ to mean that the measure $\mu$ is absolutely 
continuous with respect to the positive measure $\lambda$. 
We denoted by ${d\mu\over d\lambda}$ the 
Radon-Nikodym derivative of $\mu$ with respect to $\lambda$.

$\lp_\lambda(\Omega;\rr^N)$ stands for the usual Lebesgue space of
$p$-summable $\rr^N$-valued functions with respect to $\lambda$.
If $u\in \lu(\Omega;\rr^n)$ then $Du$ denotes its distributional gradient.
We say that $u\in\lu\orn$ is a {\it function of bounded variation}, and we 
write $u\in BV\orn$, if
all its distributional first derivatives $D_iu_j$ are Radon measure with 
finite total variation in $\Omega$;
we denote by $Du$ the $\Mn$-valued measure whose entries are $D_iu_j$.

We will use the following notion of Sobolev space with respect to a measure 
$\lambda$, which is a finite Borel positive measure on $\Om$, 
introduced by Ambrosio, Buttazzo and Fonseca \cite {ABF}
$$\wup_\lambda\orn=\Bigl\{u\in\lp\orn:\ u\in BV\orn,\ Du<\!<\lambda,\ 
{dDu\over d\lambda}\in \lp_\lambda(\Omega;\Mn)\Bigr\}$$
for all $1\le p\le +\infty$.

Let $u\in\lu\orn$, and let $Eu$ be the symmetric part of the distributional 
gradient of $u$$\ie$
$$Eu:= {E_{ij}u}, \qquad  E_{ij}u := {1 \over 2} (D_i u_j +D_j u_i).$$
The space ${\rm LD}(\Om)$ is defined as the set of all functions  
$u\in\lu\orn$ such that $ E_{ij}u\in \lu(\Om)$ for any $i,j=1,...,n$.

We say that  $u\in\lu\orn$  is a
{\it function with bounded deformation}, and we write $u\in BD(\Om)$, if 
$ E_{ij}u$ is a Radon  measure with finite total variation in $\Omega$ for any
 $i,j=1,...,n$. For every $u\in BD(\Om)$ we consider the Radon-Nikodym 
decomposition of $Eu$, with respect to the Lebesgue measure ${\cal L}^n$,
into a singular part $E^s u$ and an 
absolutely continuous part $E^a u= {\cal E}u\, {\cal L}^n$, with density
${\cal E}u = {dEu\over d{\cal L}^n}$. 
We say that $x\in \Om$ belongs to $J_u$, the jump set of $u$, if and only if 
there exist a
unit normal $\nu \in S^{n-1}$ and two vectors $a$ and $b$ in $\rr^n$
such that
$$
\lim_{\rho\to 0^+} {1\over \rho^n} \int_{B_{\rho}^+ (x,\nu)}
|u(y)-a|\dy =0
$$
$$
\lim_{\rho\to 0^+} {1\over \rho^n} \int_{B_{\rho}^- (x,\nu)} 
|u(y)-b| \dy =0
$$
where $B_{\rho}^{\pm}(x,\nu)=\{ y\in B_{\rho}(x)\ : \ 
(y-x , \pm \nu)>0 \}$ and $B_{\rho}(x)$ is the open ball of center
$x$ and radius $\rho$. The triplet $(a,b,\nu)$ is uniquely
determined up to a change of sign of $\nu$ and a permutation of
$(a,b)$. For every $x\in J_u$ we define $u^+(x)=a$, $u^-(x)=b$ and 
$\nu_u(x) =\nu$.
The singular part $E^s u$ can be written as the sum of $E^s u \res J_u$ and 
of $E^s u \res (\Om \setminus J_u)$; the first part, called the 
{\it jump part}, can be represented by 
\begin{equation}\label{decomp2}
E^s u \res J_u =(u^+ - u^-) \odot \nu_u {\cal H}^{n-1}\res J_u  
\end{equation}
while the second part, called the {\it Cantor part}, vanishes on 
any Borel set which is $\sigma$-finite with respect to ${\cal H}^{n-1}$
(see \cite{ACD} Remark 4.2 and Proposition 4.4).
We call {\it intermediate topology} on $BD(\Om)$ that defined by the
distance
\begin{equation}\label{2422}
\| u-v\|_{\lu\orn}+ ||Eu|(\Om)- |Ev|(\Om)| .
\end{equation}
If  $u\in\lu\orn$ is such that  $ E_{ij}u\in \lp(\Om)$ for any
$i,j=1,...,n$ and $\Om$ has a locally Lipschitz boundary, then we have Korn's 
inequality for all $1<p<+\infty$
\begin{equation}\label{242}
\sum_{i,j=1}^n \int_{\Om} |D_i u_j(x)|^p \dx \le
c \int_{\Om} \Bigl(|u(x)|^p + |Eu(x)|^p\Bigr)\dx\, ;
\end{equation}
hence, this space is none other than $\wup\orn$
(see Chapter 1, Section 1 in \cite{Te}).
For a general exposition of the theory of functions of bounded deformation
we refer to \cite{S}, \cite{S1}, \cite{Ko}, \cite{MSC}, \cite{AG}, 
\cite{TS},  \cite{Te}, \cite{ACD}.

If $u\in \lu\orn$, we denote by $\tilde u$ the {\it precise
representative} of $u$, whose components are defined by
\begin{equation}\label{2423}
\tilde u_i(x)=\limsup_{\rho\to 0^+}\ave_{B_{\rho}(x)}u_i(y)\dy\,.
\end{equation}

Let $f:\rr^k\to [0,+\infty]$ be a convex function. We define the
recession function $f^\infty$ of $f$ as
\begin{equation}\label{rec}
f^\infty(\xi)=\displaystyle\lim_{t\to\infty}{f(t\xi)\over t}\qquad
\mbox{for every } \xi\in\rr^k\,.
\end{equation}
It is well-known (see, for instance, \cite{Bu}) that this limit
exists, and defines a convex, subadditive and positively homogeneous
of degree one function.

We recall the definition of De Giorgi's {\it $\Gamma$-convergence} in $\lp$
spaces. Given a family of functionals
$F_j:\lp\orn\to [0,+\infty]$, $j\in\NN$, for $u\in \lp\orn$, we define
$$
\Gamma(\lp)\hbox{-}\liminf_{j\to+\infty}
F_j(u)=\inf\Bigl\{\liminf_{j\to+\infty} F_j(u_j)\ :\ u_j{
\buildrel{\lp}\over{\to}} u\Bigr\},$$
and
$$\Gamma(\lp)\hbox{-}\limsup_{j\to+\infty}
F_j(u)=\inf\Bigl\{\limsup_{j\to+\infty} F_j(u_j)\ :\ u_j{
\buildrel{\lp}\over{\to}} u\Bigr\};
$$
if these two quantities  coincide then their common value is called the
$\Gamma$-{\it limit} of the sequence $(F_j)$ at $u$,
and is denoted by
$\Gamma(\lp)\hbox{-}\lim_{j\to+\infty} F_j(u)$.
It is easy to check that $\l =\Gamma(\lp)\hbox{-}\lim_{j\to+\infty} F_j(u)$
if and only if

(a) for every sequence $(u_j)$ converging to $u$
we have
$$
\l \le\liminf_{j\to+\infty} F_j(u_j);
$$

(b) there exists a sequence $(u_j)$ converging to $u$
such that
$$
\l \ge\limsup_{j\to+\infty} F_j(u_j).
$$

We say that $(F_\e)$ $\Gamma(\lp)$-converges to $\l$  at $u$ as $\e\to 0$ 
if for
every sequence of positive numbers $(\e_j)$ converging to $0^+$ there exists a
subsequence $(\e_{j_k})$ for which we have
$$
l=\Gamma(\lp)\hbox{-}\lim_{k\to+\infty}
F_{\e_{j_k}}(u).
$$

We recall that the $\Gamma$-upper and lower limits defined above are
$\lp$-lower semicontinuous functions. For a comprehensive study of 
$\Gamma$-convergence we refer to \cite{DM} and \cite{GCB}, while a 
detailed analysis of some of its applications to homogenization theory can 
be found in \cite{BD}.

\section{The space $\ldp_{\lambda}(\Om)$}
In this section we define the analog of $\wup_{\lambda}\orn$ when the 
gradient is replaced by the linearized strain tensor.
\begin{definition}\label{24.7}
Let $\lambda$ be a finite Borel positive measure on the open set
$\Omega\subset\rn$, and let $1\le p\le+\infty$. 
We define the space

\begin{equation}\label{11}
\ldp_\lambda(\Om)=\Bigl\{u\in\lp\orn:\ u\in BD(\Om),\ Eu<\!<\lambda,\
{dEu\over d\lambda}
\in \lp_\lambda(\Omega;\Mn_{sym})\Bigr\}\,.
\end{equation}
We will use the notation ${\rm LD}_\lambda(\Om)$ instead of 
${\rm LD}^1_\lambda(\Om)$.
\end{definition}

\begin{proposition}\label{}
{\rm(i)} The spaces $\ldp_\lambda(\Om)$ and $\ldp_{\lambda'}(\Om)$ coincide 
whenever $|\lambda-\lambda'|(\Om \setminus B)=0$ for some 
${\cal H}^{n-1}$-negligible Borel subset $B$ of $\Om$.\\
{\rm(ii)} The measure $\lambda$ in Definition {\rm \ref{24.7}} can always be 
assumed
concentrated on a Borel set where its $(n-1)$-dimensional upper density is 
finite. 
\end{proposition}

{\sc Proof.} Point (i) easily follows from the fact that $BD$ functions do not 
charge
${\cal H}^{n-1}$-negligible sets (see Remark 3.3 in \cite{ACD}).
Point (ii) follows from Remark 2.3 in \cite{ABF}.
\qed

In the following proposition we prove a Leibniz-type formula for the
densities with respect to a measure $\lambda$. This formula will be used in 
the proof of the fundamental estimate, Proposition \ref{24.14}.

\begin{proposition}\label{24.8}
If $u\in \ldp_\lambda(\Om)$, $v\in \wui_\lambda(\Omega)$ and
$\tilde u \odot{dDv\over d\lambda} \in \lu_\lambda (\Om;\Mn)$ 
then $uv\in \ldp_\lambda(\Om)$, and
\begin{equation}\label{12}
{dE(uv)\over d\lambda}=
\tilde v {dEu\over d\lambda}+\tilde u \odot{dDv\over d\lambda} .
\end{equation}
\end{proposition}

{\sc Proof.} By definition, functions in $\ldp_\lambda$ have bounded
deformation.
Using the characterization of the spaces $BV(\Om)$ and $BD(\Om)$
by means of one-dimensional sections (see Proposition 3.2 in \cite{ACD})
we have
$$u^{\xi}_y \in BV(\Om^{\xi} _y), \qquad v_{y,\xi} \in BV(\Om^{\xi} _y)
\qquad  {\cal H}^{n-1}\hbox{-}a.e. \quad y\in \Om^{\xi} $$ 
where  
$$u^{\xi} _y(t) =u^{\xi}(y+t\xi)=(u(y+t\xi),\xi),\qquad 
v_{y,\xi}(t)=v(y+t\xi)
\qquad \forall t\in \Om^{\xi} _y.$$
Hence by the chain rule formula for $BV$ functions (see \cite{Bra} Section 1.8, 
\cite{AFP} Theorem 3.93 and Example 3.94) we have
$$(uv)^\xi _y=u^\xi _y\, v_{y,\xi }\in BV(\Om ^\xi _y)$$ 
and
$$D(u^{\xi} _y v_{y,\xi})=\tilde v_{y,\xi} Du^{\xi} _y+ \tilde u_y^{\xi}
Dv_{y,\xi} \qquad {\cal H}^{n-1}-a.e. \quad y\in \Om^{\xi} .$$
By Proposition 3.2 in \cite{ACD}  and by the structure theorem for $BV$
functions (see \cite{Bra} Section 1.8), we can prove that $uv\in BD(\Om)$ and 
$$ (Euv \xi, \xi) = (\tilde v Eu \xi, \xi) + (\tilde u \odot Dv \xi, \xi)
\qquad \forall \xi\in \rn.
$$
By choosing $\xi=\xi_i+\xi_j$, where $\xi_1,\dots,\xi_n$ is a basis of $\rr^n$,
we get
\begin{equation}\label{meas}
E(uv) = \tilde v Eu +\tilde u \odot Dv.
\end{equation}
Since the measures in the left hand-side of (\ref{meas}) 
are absolutely continuous with respect to $\lambda$ with densities in 
$\lp_\lambda(\Omega;\Mn_{sym})$, we finally get
 $uv\in \ldp_\lambda(\Om)$ and (\ref {12}) is proved.
\qed

\begin{remark}\label{24.09}{\rm 
Note that in (\ref{12}) it is necessary to consider the precise
representatives of $u$ and $v$, since the measure $\lambda$ may take into
account also
sets of zero Lebesgue measure.

}\end{remark}

\section{Choice of the measure and some examples}
Let $\mu$ be a non-zero  positive Radon measure on $\rn$ which
is $1$-periodic;
i.e.,
$$
\mu(B+e_i)=\mu(B)
$$
for all Borel subsets $B$ of $\rn$ and for all $i=1,\ldots,n$.
We will assume the normalization
\begin{equation}\label{2404}
\mu([0,1)^n)=1\,.
\end{equation}
For all $\e>0$ we define the $\e$-periodic positive Radon measure $\mu_\e$ by
\begin{equation}\label{2405}
\mu_\e(B)=\e^n\,\mu\Bigl({1\over\e}B\Bigr)
\end{equation}
for all Borel sets $B$. Note that by (\ref{2404}) the family $(\mu_\e)$
converges locally weakly$^*$ in the sense of measures to the Lebesgue measure
as $\e\to 0$.

In the sequel $f:\rn\times\Mn\to[0,+\infty)$ will be a fixed Borel
function
$1$-periodic in the first variable and satisfying the growth condition of order
$p\ge 1$:
there exist $0<\alpha\le\beta$ such that
\begin{equation}\label{(14.2)}
\alpha|A|^p\le f(x,A)\le \beta (1+|A|^p)\,
\end{equation}
for all $x\in\rn$ and $A\in\Mn$.

For every bounded open set $\Omega$, we define the
functionals at scale $\e>0$ as
\begin{equation}\label{2425}
F_\e(u,\Omega)=\cases{\displaystyle
\int_\Omega f\Bigl({x\over\e},{dEu\over d\mu_\e}\Bigr)d\mu_\e
& if $u\in \ldp_{\mu_\e}(\Om)$\cr
\cr
+\infty & otherwise.}
\end{equation}

Now we  consider some additional assumptions on the measure $\mu$, in order 
to prove the existence and the integral representation of the $\Gamma$-limit
of the functionals $F_{\e}$ as $\e\to 0$.
In the sequel we will point out that these conditions are necessary and 
sufficient.

We assume:

{\rm (i) (\it{existence of cut-off functions})} there exist $K>0$ and
$\delta>0$ such that for all $\e>0$, for all pairs $U,V$ of open subsets of
$\rn$ with $U\subset\subset V$, and $\dist(U,\partial V)\ge\delta\e$, and
for all $u\in \ldp_{\mu_\e}(V)$ there exists $\phi\in\wui_{\mu_\e}(V)$ with
$0\le\phi\le 1$, $\phi=1$ on $U$, $\phi=0$ in a neighbourhood of $\partial
V$, such that
\begin{equation}\label{2427}
\int_V\Bigl|{dD\phi\over d{\mu_\e}}\odot \tilde u\Bigr|^pd{\mu_\e}\le
{K\over \bigl(\dist(U,\partial V)\bigr)^p}\int_{V\setminus U}|u|^p\dx\,.
\end{equation}
Such a $\phi$ will be called a {\it cut-off function between} $U$ and $V$;

{\rm (ii) ({\it existence of periodic test-functions})} for all 
 $i,j=1,\ldots,n$, there
exists $z_{ij}\in \ldp_{\mu,\loc}(\rn)$ such that
 $x\mapsto z_{ij}(x)-x_j e_i$ is
$1$-periodic.

\begin{remark}\label {24.11a}{\rm
Note that if $\mu$ is $p$-homogenizable in the sense of Definition 3.2 
in \cite{ABC}; {\it i.e., }if there exists $z_i \in \wup_{\mu,\loc} (\rn)$
such that  $x\mapsto z_{i}(x)-x_i$ is $1$-periodic, 
then the functions $z_{ij}=z_j e_i$ trivially satisfy the condition
(ii) above but the converse is not true.
}\end{remark}

\begin{remark}\label{24.11}{\rm
Note that the Lebesgue measure trivially satisfies properties 
(i), (ii). Note that property (i)  depends on $\mu$ and $p$.
}\end{remark}

We consider in our context the measure $\mu$ of Examples 3.1(a) and (b) in
\cite{ABC}. 

\begin{example}\label{24.9}{\rm ({\it
Perfectly-rigid bodies connected by springs.})\\
We consider
$$
E=\{y\in\rn:\ \exists i\in\{1,\ldots, n\}\hbox{ such that } \ y_i\in\ZZ\},
$$
that is, the union of all the boundaries of cubes $Q_i=i+(0,1)^n$ with
$i\in\Zn$.
$E$ is an $(n-1)$-dimensional set in $\rn$. We set
$$
\mu(B)={1\over n}{\cal H}^{n-1}(B\cap E)
$$
for all Borel sets $B$.
For every $\e>0$ we have
$$
\mu_\e(B)={1\over n} \e{\cal H}^{n-1}(B\cap \e E)
\,.
$$
If $u\in \ldp_{\mu_\e}(\Om)$ then  
$Eu=0$ on every connected component of each  $\e Q_i\cap\Omega$, so 
in this case $\ldp_{\mu_\e}(\Om)$ consists of  functions which are rigid 
displacements 
on these sets\ie $u_i=R_i x + c_i$  on each $\e Q_i\cap\Omega$ with $R_i$
a  $n\times n$ skew symmetric matrix, and $c_i \in \rr^n$. Hence by 
(\ref{decomp2}), we have
$$
{dEu\over d\mu_\e}={n\over \e}{d Eu\over d{\cal H}^{n-1}}
={n\over\e}(u_i-u_j)\odot(i-j) \hbox{ on }
\partial(\e Q_i)\cap \partial(\e Q_j)\cap \Omega\,.
$$
In this case the functionals
$F_\e$ take the form
$$
\e\int_{\Omega\cap\e E}g\Bigl({x\over\e},{1\over\e}{d Eu\over d{\cal H}^{n-1}}
\Bigr)
d{\cal H}^{n-1}.                 
$$
Note that if
$\Omega$ is bounded then $\ldp_{\mu_\e}(\Om)={\rm LD}^{\infty}_{\mu_\e}(\Om)$
for all
$p$ if the number of connected components of each $\Omega\cap\e Q_i$ is 
finite.

Comparing with Example 3.1(a) in \cite{ABC}, we get that
 $ \wup_{\mu_\e}\orn \subset \ldp_{\mu_\e}(\Om) $.

The measure $\mu$ satisfies the conditions (i) and (ii)
  for all $p\ge 1$.  In fact, to prove (i) we consider the same cut-off
function in Example 3.4(a) of \cite {ABC}
$$
\phi(x)=1-\Bigl({1\over C}\,\Bigl[{1\over\e}\inf\{|x-y|_\infty: y\in
U_\e\}\Bigr] \wedge 1\Bigr),
$$
where fixed $\e>0$, $U_\e=\bigcup\{\e Q_i:\ \e Q_i\cap U\neq\emptyset\}$,
$|x-y|_\infty=\max_{1\le i\le n}|x_i-y_i|$, 
and
$$
C=\Bigl[{1\over\e}\inf\Bigl\{|x-y|_\infty: x\in U_\e,\ y\in\partial
V\Bigr\}\Bigr]-2
$$
(we denote $[t]$ the integer part of $t$).
Note that $|dD\phi/d\mu_\e|\le n/(C\e)\le c/\dist(U, \partial V)$ for some
constant $c$ independent of $U$ and $V$.

Interpreting $u^{\pm}$ as traces of Sobolev functions defined on
each cube $Q_i$, we have
$$
\left(\int_{\partial Q_i}|u^{\pm}|^p
d{\cal H}^{n-1}\right)^{1/p} \leq c
\|u\|_{\wup(Q_i)},
$$
hence by a scaling argument and by Korn's inequality (\ref{242})
\begin{eqnarray*}
\left(\e \int_{\partial\e Q_i}|u^{\pm}|^p
d{\cal H}^{n-1}\right)^{1/p} 
&\leq& 
c \left(\int_{\e Q_i}|u|^p
dx\right)^{1/p} +\left({1 \over \e}\int_{\e Q_i}|Eu|^p
dx\right)^{1/p} \\
&=& 
c \left(\int_{\e Q_i}|u|^p
dx\right)^{1/p}
\end{eqnarray*}
where $c$ depends only on the cube.
If $p=1$ we can apply the trace inequality in $\ld(Q_i)$ 
$$
 \int_{\partial Q_i}|u^{\pm}|
d{\cal H}^{n-1} \leq c
\int_{Q_i}|u|
dx +|Eu|(Q_i)\, ,$$
 so we get

$$\e \int_{\partial\e Q_i}| u^{\pm}|
d{\cal H}^{n-1} \leq c
\int_{\e Q_i}|u|dx.$$
Hence for all $p\ge 1$
$$\e \int_{\partial\e Q_i}| u^{\pm}|^p\,
d{\cal H}^{n-1} \le c
\int_{\e Q_i}|u|^p dx.$$
For two cubes

$$
\e \int_{\partial \e Q_i\cap \partial \e Q_j}|\widetilde u|^p\,d{\cal
H}^{n-1} \le
\e \int_{\partial \e Q_i\cap \partial \e Q_j}(|u_i|^p+|u_j|^p)\,d{\cal
H}^{n-1} \le c \int_{ \e Q_i\cup \e Q_j}|u|^p\dx\,
$$
so that
\begin{eqnarray*}
\int_V\Bigl|{dD\phi\over d\mu_\e}\odot \widetilde u\Bigr|^p d\mu_\e
&\le&
{c^p\e\over \dist(U, \partial V)^p} \int_{(V\setminus U)\cap \e E\cap
\spt\!D\phi}|\widetilde u|^p \,d{\cal H}^{n-1}\\
&\le &
2n{c^p\over \dist(U, \partial V)^p}\int_{V\setminus U}|u|^p\dx\,.
\end{eqnarray*}
The proof of (i) is then complete.
To verify (ii) we apply  Remark~\ref{24.11a} to Example 3.4(a)
in \cite {ABC} and take simply $z_{ij}(x)= [x_j]e_i$.
}\end{example}

\begin{example}\label{24.9bis}{\rm ({\it 
Elastic media connected by springs}).\\
Let $E$ be as in the previous example and let
\begin{eqnarray*}
\mu(B)&=&{1\over n+1} \Bigl(|B|+ {\cal H}^{n-1}(E\cap B)\Bigr)
\\
\mu_\e(B)&=&{1\over n+1}\Bigl(|B|+ \e{\cal H}^{n-1}((\e E)\cap B)\Bigr)\,.
\end{eqnarray*}
In this case the functions in $\ldp_{\mu_\e}(\Om)$ are functions whose
restriction to each $\e Q_i\cap\Omega$ belongs to $\wup(\e
Q_i\cap\Omega;\rn)$ 
when $p>1$ by the Korn's inequality (\ref{242}) (we suppose that 
$\e Q_i\cap \Om$ has a
locally Lipschitz boundary) and to $\ld(\e Q_i\cap\Omega)$ when $p=1$, while 
the difference of the traces on both sides of
$\partial(\e Q_i)\cap \partial(\e Q_j)\cap \Omega$ is $p$-summable for every
$i,j\in\Zn$. Hence if we compare our case with Example 3.1(b) in \cite{ABC},  
we can conclude that $ \wup_{\mu_\e}\orn =\ldp_{\mu_\e}(\Om) $
if $p>1$
 and $\wuu_{\mu_\e}\orn \subset \ld_{\mu_\e}(\Omega)$ if $p=1$.
The functionals $F_\e$
take the form
$$
{1\over n+1}\int_{\Omega}f\Bigl({x\over\e},{d Eu\over dx}\Bigr)\dx+
\e\int_{\Omega\cap\e E}g\Bigl({x\over\e},{1\over\e}{d Eu\over d{\cal
H}^{n-1}}\Bigr)d{\cal
H}^{n-1}.
$$
The measure $\mu$  satisfies conditions (i) and (ii)  for
all $p\ge 1$ by  Example 3.4(b) in \cite {ABC}.

}\end{example}

\section {The homogenization theorem}
The homogenization theorem for the functionals in (\ref{2425})
takes the following form.

\begin{theorem}\label{24.12}
Let $\mu$ be a measure which satisfies conditions {\rm(i)} and {\rm(ii)}
in Section {\rm 4}, 
and for every bounded open subset
$\Omega$ of $\rn$ let $F_\e(\cdot,\Omega)$ be defined on $\lp\orn$ by
{\rm(\ref{2425})}. Then the $\Gamma$-limit
\begin{equation}\label{2428}
F_\Hom(u,\Omega)=\Gamma(\lp)\hbox{-}\lim_{\e\to 0}F_\e(u,\Omega)
\end{equation}
exists for all bounded open subsets $\Omega$ with Lipschitz boundary
and for all $u\in\lp\orn$; it
can be represented on $\wup\orn$ for $p\ge 1$ as
\begin{equation}\label{2429} F_\Hom(u,\Omega)=\int_\Omega
f_\Hom(Eu)\dx\,, \end{equation}
where the homogenized integrand satisfies
the asymptotic formula
\begin{eqnarray}\label{2430}
f_\Hom(A)&=&
\lim_{k\to+\infty} \inf\Bigl\{{1\over k^n}\int_{[0,k)^n} f\Bigl(x,{dEu\over
d\mu}\Bigr)d\mu:
\\ \nonumber&&\qquad\qquad
u\in \ldp_{\mu,\loc}(\rn),\
u- A x\ \
\hbox{$k$-periodic}\Bigr\}\,
\end{eqnarray}
for all $A\in \Mn_{sym}$.

Moreover, $F_\Hom(u,\Omega)=+\infty$ if $p>1$ and 
$u\in\lp\orn\setminus\wup\orn$, or if $u\in \lu\orn \setminus BD(\Om)$ when 
$p=1$.

Furthermore, if $f$ is convex then the  $\Gamma$-limit can be
 represented as
$$F_\Hom(u,\Omega)=\int_\Omega
f_\Hom({\cal E} u)\dx +\int_\Omega
f_\Hom^{\infty}\Bigl({dEu^s\over d|Eu^s|}\Bigr)d|Eu^s|$$
for all $u\in BD(\Om)$ when $p=1$.
\end{theorem}

\begin{remark}\label{24.12bis}{\rm
Note that we cannot replace the sets $[0,k)^n$ by the sets 
$(0,k)^n$ 
if $\mu([0,k)^n\setminus (0,k)^n) \neq 0$, see  Remark 3.6 in \cite{ABC}.

Same examples and considerations of Remarks 3.7 and 3.8 in \cite{ABC}, 
applied to our case, show that condition (ii) for the measure $\mu$ and
the assumption that $\Om$ has a Lipschitz boundary are necessary to get a 
homogenization theorem.  
In fact, if condition (ii) fails then $f_\Hom(A)=+\infty$ if $A\neq 0$; while 
if $\Omega$ does not have Lipschitz boundary then the equality (\ref{2429}) 
may not hold.
}\end{remark}

The following proposition is a usual tool to prove the existence of the 
$\Gamma$-limit and its integral representation (see \cite{DM} Chapter 18, 
\cite{BD} Chapter 11).
\begin{proposition}[Fundamental Estimate]\label{24.14}
For every $\sigma>0$ there exists $\e_\sigma$ and $M>0$
such that for all $U, U', V$ open subsets of $\Omega$ with
$U'\subset U$ and $\dist(U',V\setminus U)>0$,
for all $\e<\e_\sigma\dist(U',V\setminus U)$
and for all $u\in \ldp_{\mu_\e}(\Om)$,
$v\in \ldp_{\mu_\e}(\Om)$ there exists a cut-off function between
$U'$ and
$U$, $\phi\in\wui_{\mu_\e}(U\cup V)$, such that
\begin{eqnarray}\label{2433}
F_\e(\phi u+(1-\phi)v,U'\cup V)&\le&
(1+\sigma)(F_\e(u,U)+F_\e(v,V))
\\ \nonumber
&&\hskip-3cm
+{M\over\bigl(\dist(U',V\setminus U)\bigr)^p}\int_{(U\cap V)\setminus U'}
|u-v|^pdx+\sigma{\mu_\e}((U\cap V)\setminus U').
\end{eqnarray}
\end{proposition}

{\sc Proof.} By taking (\ref{12}) and condition (i) into account, the proof 
follows
exactly that of Proposition 4.1 \cite{ABC}.

\begin{proposition}\label{24.13}
For every $A\in\Mn_{sym}$ there exists $z_A\in \ldp_{\mu,\loc}(\rn)$ 
such that
$z_A- A x$ is $1$-periodic and satisfies
\begin{equation}\label{2432}
\int_{[0,1)^n}\Bigl|{dEz_A\over d\mu}\Bigr|^pd\mu\le
 c |A|^p.
\end{equation}
\end{proposition}

{\sc Proof.} Define $z_A=\sum_{i,j=1}^n A_{ij} z_{ij} $, 
where $z_{ij}$
are as in condition (ii). Inequality (\ref{2432}) is then trivial.
\qed

We fix $(\e_j)$ which goes to zero. We define
$$
F'(u,U)=\Gamma(\lp)\hbox{-}\liminf_{j\to+\infty} F_{\e_j}(u,U)
$$
$$
F''(u,U)=\Gamma(\lp)\hbox{-}\limsup_{j\to+\infty} F_{\e_j}(u,U)
$$
for all $u\in \lp\orn$ and for all open subsets $U$ of $\Omega$.

\begin{proposition}[Growth Condition]\label{24.14bis}
We have for all open subsets $U$ of $\Omega$
with $|\partial U|=0$

$$
F''(u,U)\le c \int_U(1+|Eu|^p)dx
$$
for all $u\in\wup\orn$ if $p>1$ and  

$$
F''(u,U)\le c(|U|+|Eu|(U))
$$
for all $u\in BD(\Om)$ if $p=1$.
\end{proposition}

{\sc Proof.} This Growth Conditions can be obtained
modifying the proof of Proposition 4.3 in \cite {ABC}.
In particular in Step 2 therein now we have to consider the  affine 
functions
$u_i(x)=A_i x+c_i$ for some $A_i\in\Mn_{sym}$
and $c_i\in\rn$, in Step 3 we just have to note that piecewise
affine functions are dense in $BD$ endowed with the intermediate
topology (\ref{2422})
(see \cite{Te} Theorem 3.2 Chapter $2$ Section $3$).
\qed

\begin{proposition}\label{24.15}
There exists a subsequence of $(\e_j)$ (not relabeled) such that
for all open subsets $U$ of $\Omega$  with $|\partial U|=0$
 there exists the $\Gamma$-limit
$$
\Gamma\hbox{-}\lim_{j\to+\infty} F_{\e_j}(u,U)=F(u,U)\,,
$$
for all $u\in\wup\orn$ if $p>1$ and for all $u\in BD(\Omega)$ if $p=1$.
There exists a function $\varphi:\Mn\to\rr$ such that
$$
F(u,U)=\int_U\varphi(Eu)dx
$$
for all $u\in\wup\orn$ if $p\ge 1$; moreover if f is convex 
$$
F(u,U)=\int_U\varphi\Bigl({\cal E}u\Bigr)dx+\int_U\varphi^{\infty}
\Bigl({dE^s u\over d|E^s u|}\Bigr)d|E^s u|
$$
for all $u\in BD(\Omega)$ if $p=1$.
\end{proposition}

{\sc Proof.} To prove the existence of the $\Gamma$-limit on $\wup\orn$
for $p>1$ and $BD(\Om)$ for $p=1$, 
and the integral representation of the $\Gamma$-limit
$$
F(u,U)=\int_U\varphi(Du)dx
$$
on $\wup\orn$ when $p\ge1$, we repeat the proof of Proposition 4.4
\cite{ABC} using 
Propositions \ref{24.14} and \ref{24.14bis}.
Moreover, we can prove that $\varphi (Du)= \varphi (Eu)$.
In fact, let $w_j \to Ax$ be such that 
$$F(Ax,\Om)=\lim_{j \to +\infty} F_{\e_j}(w_j,\Om)$$
and let $Rx+c$ be a rigid displacement, then
\begin{eqnarray*}
F(Ax+Rx+c,\Om) &\le& \liminf_{j\to+\infty}F_{\e_j}(w_j+Rx+c,\Om)\\
&=&
\lim_{j\to+\infty}F_{\e_j}(w_j,\Om)=F(Ax,\Om)
\end{eqnarray*}
so that $\varphi(A+R)\le \varphi(A)$.
The reverse inequality follows similarly, therefore for all $R$ $(n\times n)$ 
skew-symmetric
matrix 
$$\varphi(A+R)=\varphi(A)$$
which implies $\varphi(B)=\varphi({B+B^T\over 2})$ for any $B\in \Mn$.

Let us prove the integral representation of the $\Gamma$-limit
on $BD(\Om)$ whenever $f$ is convex.
We consider the functional defined in ${\rm L}^1_{\loc}\orn$ 
$$
G(u)=\cases{\displaystyle
\int_\Omega \varphi(Eu)\dx
& if $u\in {\rm C}^1\orn$\cr
\cr
+\infty & otherwise,}
$$
and we introduce 
$$
\overline{G}(u,U)=\inf \Bigl\{\liminf_{h\to +\infty} G(u_h, U) : u_h\in 
{\rm C}^1\orn \quad u_h \rightarrow u \quad\hbox{in}\quad {\rm L}^1_{\loc}
\orn\Bigr\}
$$
the relaxed functional of $G$. 
It is well known that $\varphi$ is convex and it is easy to check that
$\varphi(A)\ge c |A|$ for every $A\in \Mn_{sym}$, hence by the lower
semicontinuity and relaxation theorems for functionals of measures 
(see for instance \cite{GS}, \cite{Bu}), we obtain
$$
\overline{G}(u,U) = 
\int_U\varphi({\cal E}u)\dx+\int_U\varphi^{\infty}
\Bigl({dE^s u\over d|E^s u|}\Bigr)d|E^s u|
$$
for every $u\in BD(\Om)$ (see \cite{Te} Section 5).
Since $F(\cdot, U) \le G(\cdot,U)$ in ${\rm W}^{1,1}\orn$, 
by the lower semicontinuity of the $\Gamma$-limit we obtain
$$
F(u,U)\le \int_U\varphi({\cal E}u)\dx+\int_U\varphi^{\infty}
\Bigl({dE^s u\over d|E^s u|}\Bigr)d|E^s u|
$$
for all $u\in BD(\Om)$.
The reverse inequality is obtained by a convolution argument. In fact 
we consider 
$U_k =\{ x\in U : d(x,\partial U)> {1\over k}\}$, $\rho_k$ with 
$\spt\rho _k \subset B(0,{1\over k})$ and $u_k=u*\rho_k$.
For $y\in  B(0,{1\over k})$ and $k$ large enough we have
that $U_k \subset y+U$.\\ 
Since $F(\cdot, U)$ is convex for all  $U\in\Ao$ 
and $F(u^y, U_k) \le F(u,U)$ with $u^y(x)=u(x-y)$, by Jensen's inequality 
$$
F(u*\rho_k, U_k)\le F(u,U).
$$
On the other hand, we also have
$$
\lim_{k\to +\infty} F(u_k,U_k)= \overline{G}(u,U)
$$
hence  we can conclude that 
$$
F(u,U)=\int_U\varphi({\cal E}u)\dx+\int_U\varphi^{\infty}
\Bigl({dE^s u\over d|E^s u|}\Bigr)d|E^s u|
$$
as desired.
\qed

\begin{proposition}[Homogenization Formula]\label{24.16}
For all $A\in\Mn_{sym}$ there exists the limit in {\rm(\ref{2430})}
and we have $\varphi(A)=f_\Hom(A)$.
\end{proposition}

{\sc Proof.} It can be obtain repeating the proof of the Proposition 4.5 of 
\cite{ABC}
but defining 
$$
g_k(A)=\inf\Bigl\{{1\over k^n}\int_{{(0,k)^n}}
f(x,{dEu\over d\mu})d\mu:
u\in \ldp_{\mu,\loc}(\rn),\
u-A  x\ \
\hbox{$k$-periodic}\Bigr\}
$$ 
for all $A\in \Mn_{sym}$ and $k\in \NN$.
\qed

{\sc Proof of Theorem \ref{24.12}.} It remains to check the coercivity of the 
$\Gamma$-limit.
By the growth condition on $f$ and a comparison argument, it is enough
to prove this for $f(A)=|A|^p$. 
We know that  the
$\Gamma$-limit $F_\Hom$ exists for all $u\in \lp\orn$ and for all sets 
$R$ in the countable family
${\cal R}$ of all finite unions of open rectangles of $\Omega$ with rational
vertices, in this case $F_\Hom$ is also convex.
For all $U'$, $U\in\Ao$ such that $U'\subset\subset U$ there exists
$R\in \cal R$ such that $U'\subset\subset R \subset\subset U$.
Reasoning as in the previous proof, for $y\in  B(0,{1\over k})$ and $k$ 
large enough we have that $R \subset y+U$ hence 
$$F_\Hom(u_k, R)\le F'(u,U)$$
and 
\begin{equation}\label{ine}
\liminf_{k\to +\infty} F_\Hom(u_k, U')\le F'(u,U)
\end{equation}
with $u_k=u*\rho_k$ (see \cite{DM} Chapter 23).

It will be enough then to prove that $f_\Hom(A)\ge c|A|^p$.
In fact for any $u\in\lp\orn\setminus\wup\orn$ when $p>1$ by (\ref{ine}) 
$$F'(u,U)\ge c \liminf_{k\to +\infty} \int_{U'} |Du_k|^p \dx$$
by the arbitrarity of $U'$, we get $F_\Hom(u,U)=+\infty$.
Similarly, if $p=1$ for all $u\in \lu\orn \setminus BD(\Om)$
we have $|Eu|(\Omega)=+\infty$, let $\Omega '\subset \subset \Omega$ 
we get by (\ref{ine}) that
$$F'(u,\Omega)\ge c \liminf_{k\to +\infty}|Eu_k|(\Omega')$$
by arbitrarity of $\Omega'$ we obtain $F_\Hom (u,\Omega) =+\infty$.

Since $f_\Hom$ is positively homogeneous of
degree $p$, to prove that $f_\Hom(A)\ge c|A|^p$, it is sufficient 
to check that $f_\Hom(A)\neq 0$ if $A\neq 0$.
To this
aim, let $u_\e\to Ax$ be such that $F_\e(u_\e,(0,1)^n)\to f_\Hom(A)$.
If $f_\Hom(A)=0$ then by a ``Poincar\'e-type" inequality for $BD$ functions
(Proposition 2.3 Chapter 2 of \cite{Te}), by H\"older's inequality and a
scaling argument we obtain that 
\begin{eqnarray*}
0= f_\Hom(A) &=& \lim_{\e \to 0} \int_{(0,1)^n}
 \Bigl |{dEu_{\e} \over d\mu_{\e}}\Bigr|^p\, d\mu_{\e}
\\ 
&\ge&
\lim_{\e\to 0} c \Bigl(\int_{(0,1)^n}|u_{\e}-Ru_{\e}| \dx\Bigr)^p
\end{eqnarray*}
where the constant $c$ depends only on $\Om$ and $Ru_{\e}$ is a rigid 
displacement. Hence 
$Ru_{\e} \to Ax$ in $\lu$, and we get 
 a contradiction because $A$ is a symmetric matrix.
\qed

\section{Non local effects}
Theorem \ref{24.12} shows the  $\Gamma(\lp)$-convergence of the  
functionals $F_{\e}$ to $F_{\hom}$ in $\wup\orn$ and that the $\Gamma$-limit 
is local; in fact 
we have represented $F_{\hom}$ as the integration over $\Om$ of a local 
density of energy of the form $f_{\hom}(Eu)$.

Now, if we consider
$$
F^\gamma_\e(u,\Om)= \e^{\gamma} \int_{\Om} f\Bigl({dEu\over d\mu_\e}\Bigr) 
d\mu_\e
$$
then $\Gamma(\lp)\hbox{-}\lim_{\e\to0}F^\gamma_\e(u,\Om)=0$ on $\wup\orn$,
when $\gamma>0$.
In this case, however, no coerciveness result may hold for sequences $(u_{\e})$
with $\sup_{\e>0} F^\gamma_\e(u_{\e},\Om)< +\infty$ in any norm.

We will show with an example that a more complex notion of convergence may have 
to be introduced and that
the $\Gamma$-limit functionals may be of a non-local nature.

Let $\Om= \omega \times (0,1)$ be a `cylindrical' domain where 
$\omega$ is a connected open subset of $\rr^2$.

We define
$\e D_i$ to be a two dimensional disk centered at 
$x_i= (\e i_1+ {\e\over 2}, \e i_2+ {\e\over 2})$ of radius $\e /4$ 
$$\e {E_i}^2= \e D_i \times (0,1)\qquad \qquad 
\e E^2= \bigcup_{i\in I_{\e}}\e {E_i}^2
$$
where $i=(i_1,i_2)\in I_{\e}=\{i\in \ZZ^2 : \e {E_i}^2 \subset \Om\}$, 
$$\e E^1= \Om\setminus \e E^2.$$
We call $E= D_0 \times (0,1)$.

We consider the measures
$$
\mu_\e(B)= \e{\cal H}^2 (B\cap \partial \e E^2)\, 
$$
and the functionals
$$
F^\gamma_\e(u,\Om)=\e^{\gamma}\int_{\Omega} \Bigl|{d Eu\over d\mu_\e}
\Bigr|^2 d\mu_\e \,. 
$$
Note that, up to normalization, $\mu_{\e}$ is the same measure of Example 
\ref{24.9}.

In this case $\ld^2_{\mu_\e}(\Om)$ consists of  functions which are
rigid displacements on the sets $\e E^1$ and $ \e E^2$\ie$u\in \ld^2_{\mu_\e}
(\Om)$
if and only if there exist $a_i , b_i, c, d \in \rr^3$ such that

$$u=c \wedge x +d  \quad\hbox{ on } \quad\e E^1 $$
$$u=a_i \wedge x +b_i  \quad\hbox{ on } \quad\e {E_i}^2$$
for each $i\in I_{\e}$. 
We use the notation $x=(x_{\alpha},x_3) \in\rr^3$, $x_{\alpha}=(x_1,x_2)$.

Hence 
$$
{dEu\over d\mu_\e}={1\over \e}{d Eu\over d{\cal H}^2}
={1\over\e}( c \wedge x +d - a_i \wedge x - b_i)\odot \nu \hbox{ on }
\partial(\e E_i^2)\,.
$$

\begin{definition}\label{def}

Let $u_{\e}\in \ld^2_{\mu_\e}(\Om)$. We say that $u_{\e}$ converges to 
$(u_1 , u_2)\in {\rm L}^2(\Om;\rr^3)\times {\rm L}^2(\Om;\rr^3)$ if and only 
if  
\begin{eqnarray}\label{b2}
\lim_{\e \to 0} \int_{\e E^1} |u_{\e} -u_1|^2 \dx\, =0
\end{eqnarray}
\begin{eqnarray}\label{b1}
\lim_{\e \to 0} \int_{\e E^2} |u_{\e} -u_2|^2 \dx\, =0\, .
\end{eqnarray}
\end{definition}

We will study the $\Gamma$-limit $F$ of $F^\gamma_\e$ with respect to the
convergence introduced in 
Definition \ref{def} (see Theorem \ref{teorema}). The domain of $F$ 
will be the set of pairs $(u_1,u_2)$ such that $u_1$ is a rigid displacement
and $u_2$ is in the space ${\cal U}$ of functions whose `vertical sections
are rigid displacements', introduced in the following proposition.
 
Let us define, for all $\eta>0$, 
$T_{\eta}^k = Q_{\eta}^k\times (0,1)$ where $Q_{\eta}^k= k + (0,\eta)^2$
with 
$k=(k_1,k_2)\in J=\{k\in \ZZ^2 : T_{\eta}^k\cap \Omega \neq \emptyset\}$.

\begin{proposition}
Let $u_{\e}\in \ld^2_{\mu_\e}(\Om)$ and $u_2\in {\rm L}^2(\Om;\rr^3)$.
$$
\lim_{\e \to 0} \int_{\e E^2} |u_{\e} -u_2|^2 \dx\, =0
$$
if and only if $u_2\in {\cal U}$ where
\begin{eqnarray*}
&{\cal U}&=
\Bigl\{v\in {\rm L}^2(\Omega;\rr^3) : \forall \eta>0\; \exists J\subset \ZZ^2
\; \hbox{and}\; \exists \, A^k\wedge x+ B^k\;
\hbox{on}\; T_{\eta}^k\quad \forall k\in J \; \hbox{such that} \\
\nonumber&& \qquad\qquad\qquad
\bigcup_{k\in J}T^k_{\eta}\cap \Omega = \Omega \quad\hbox{and}\quad 
\sum_{k\in J}\int_{T_{\eta}^k\cap \Om}
|v(x) - A^k\wedge x -B^k|^2\dx \le o(\eta)
\Bigr\}\, .
\end{eqnarray*}
\end{proposition}

{\sc Proof.} Let $u_{\e}\in \ld^2_{\mu_\e}(\Om)$, by definition $u_{\e}=
a_{\e,i} \wedge x + b_{\e,i} $ on $\e E^2_i$. 
Let $h\in \NN$ and $\eta >0$ such that $\eta = h\e$, we extend 
$a_{\e,i} \wedge x + b_{\e,i}$ to 
$T_{\eta}^k$ for each $i\in I_k = \{i\in \ZZ^2 : \e E^2_i \subset 
T_{\eta}^k \}$, hence we can construct a rigid displacement 
on $T_{\eta}^k$
$$
A_{\e}^k\wedge x +B_{\e}^k = {1\over h^2}\sum_{i\in I_k} a_{\e,i} \wedge x +
b_{\e,i}\,.
$$

Let us suppose that $u_{\e}$ satisfies condition (\ref{b1}),
\begin{eqnarray}\label{gamma}
&&
\int_{T_{\eta}^k \cap \e E^2} \Bigl| u_2(x) - A_{\e}^k\wedge x -
B_{\e}^k \Bigr|^2
\dx \\ \nonumber
&\le& 
c \Bigl(\sum_{j\in I_k}\int_{\e E^2_j} \Bigl| u_2(x) - a_{\e,j} \wedge x -
b_{\e,j}\Bigr|^2
\dx \\ \nonumber
&&
+\sum_{j\in I_k}\int_{\e E^2_j} \Bigl| a_{\e,j} \wedge x +b_{\e,j}-
{1\over h^2}\sum_{i\in I_k} a_{\e,i} \wedge x +
b_{\e,i}\Bigr|^2\dx \Bigr)\,.
\end{eqnarray}
Let us estimate the last term in (\ref{gamma})
\begin{eqnarray*}
&&
\sum_{j\in I_k}\int_{\e E^2_j} \Bigl| a_{\e,j} \wedge x +b_{\e,j}-
{1\over h^2}\sum_{i\in I_k} a_{\e,i} \wedge x +
b_{\e,i}\Bigr|^2\dx \\
&\le&
c \Bigl(\sum_{j\in I_k}\int_{\e E^2_j} \Bigl|a_{\e,j} \wedge x +
b_{\e,j} - u_2(x)\Bigr|^2\dx \\
&&
+\sum_{j\in I_k}\int_{\e E^2_j} \Bigl|
{1\over h^2}\sum_{i\in I_k} a_{\e,i} \wedge x + b_{\e,i} - 
u_2(x+ x_i-x_j)\Bigr|^2\dx \\
&&
+\sum_{j\in I_k}\int_{\e E^2_j} \Bigl| 
{1\over h^2}\sum_{i\in I_k} u_2(x) - u_2(x+x_i-x_j)\Bigr|^2\dx \Bigr)\,.
\end{eqnarray*}
For each $x\in \e E^2_j$ we have that $x+x_i-x_j\in \e E^2_i$, hence with 
a change of coordinates we get
\begin{eqnarray}\label{beta}
\nonumber&&
\sum_{j\in I_k}\int_{\e E^2_j} \Bigl| a_{\e,j} \wedge x +b_{\e,j}-
{1\over h^2}\sum_{i\in I_k} a_{\e,i} \wedge x +
b_{\e,i}\Bigr|^2\dx \\ \nonumber
&\le&
c \Bigl(\sum_{j\in I_k}\int_{\e E^2_j} \Bigl|a_{\e,j} \wedge x +b_{\e,j}-u_2(x)
\Bigr|^2 
\dx\\ \nonumber
&&
+\sum_{i,j\in I_k}{1\over h^2}\int_{\e E^2_i} \Bigl| 
a_{\e,i} \wedge (x+x_j-x_i) +b_{\e,i}- u_2(x)\Bigr|^2 \dx \\ \nonumber
&&
+\sum_{i,j\in I_k}{1\over h^2}\int_{\e E^2_j} \Bigl| u_2(x)- u_2(x+x_i-x_j)
\Bigr|^2 \dx \Bigr)\\ 
&\le&
c \Bigl(\sum_{i\in I_k}\int_{\e E^2_i} \Bigl|a_{\e,i} \wedge x +b_{\e,i}- 
u_2\Bigr|^2
\dx \\ \nonumber
&&
+\sum_{i,j\in I_k}{1\over h^2}\int_{\e E^2_i} \Bigl| 
a_{\e,i} \wedge (x_j-x_i)\Bigr|^2 \dx \\ \nonumber
&&
+\sum_{i,j\in I_k}{1\over h^2}\int_{\e E^2_j} \Bigl| u_2(x)- u_2(x+x_i-x_j)
\Bigr|^2 \dx \Bigr)\,.
\end{eqnarray}

Now if we denote $\Lambda$ the set of all translations of the type $x_i-x_j$ 
with $i,j\in I_k$ we get that
\begin{eqnarray}\label{alpha}
&&
\sum_{i,j\in I_k}{1\over h^2}\int_{\e E^2_j} \Bigl| u_2(x)- u_2(x+x_i-x_j)
\Bigr|^2 \dx \\ \nonumber
&\le&
\sum_{\tau\in \Lambda} {1\over h^2}\sum_{r\in C(k)}\int_{T_{\eta}^r}
|u_2(x)-u_2(x+\tau)|^2\dx 
\end{eqnarray}
where 
$C(k)=\{(k_1,k_2), (k_1\pm 1,k_2), (k_1,k_2\pm 1), (k_1\pm 1,k_2\pm 1)\}$.

Since $|\Lambda|= c\, h^2$, by (\ref{alpha}) we have
\begin{eqnarray}\label{trasl}
\nonumber&&
\sum_{k\in J}\sum_{i,j\in I_k}{1\over h^2}\int_{\e E^2_j} 
\Bigl|u_2(x)- u_2(x+x_i-x_j)\Bigr|^2 \dx \\ \nonumber
&\le&
c \sum_{\tau\in \Lambda} {1\over h^2} \Vert u_2(\cdot)-u_2(\cdot +\tau)
\Vert^2_{{\rm L}^2(\Omega;\rr^3)}\\
&\le&
c \sup_{|\tau|\le \sqrt{2}\eta}\Vert u_2(\cdot)-u_2(\cdot +\tau)
\Vert^2_{{\rm L}^2(\Omega;\rr^3)}\,.
\end{eqnarray}

Let us consider the cubes $Q_{\e,i}^j=(\e i + (0, 1)^2 )
\times (\e j + (0,\e))$ for $i\in I_{\e}$, and
$j\in J_{\e}=\{j\in \ZZ : Q_{\e,i}^j \cap \e E^2_i \neq \emptyset\}$.
Since $u_2\in {\rm L}^2(\Om;\rr^3)$, we can assume that there exists
a sequence $(u_{\e,2})$  
which is constant on each $Q_{\e,i}^j$ such that
\begin{eqnarray}\label{b22}
\lim_{\e \to 0} \int_{\Om} |u_2 - u_{\e,2}|^2 \dx =
\lim_{\e \to 0} \sum_{i\in I_{\e}}\sum_{j\in J_{\e}}\int_{Q_{\e,i}^j \cap \Om} 
|u_2 - u_{\e,2,i,j}|^2 \dx = 0
\end{eqnarray}
where $ u_{\e,2,i,j}$ is the value of $(u_{\e,2})$ on $Q_{\e,i}^j$.

So by (\ref{b1}) ) we get  
\begin{eqnarray}\label{a}
\lim_{\e \to 0} \sum_{i\in I_{\e}}\sum_{j\in J_{\e}} \int_{Q_{\e,i}^j \cap 
\e E_i^2} 
|u_{\e} - u_{\e,2,i,j}|^2 \dx = 0 \, .
\end{eqnarray}

Note that the ${\rm L}^2$-norm on the set ${\cal R}$ of rigid displacements is 
equivalent to the norm on ${\cal R}$ 
$$
\Vert a \wedge x + b\Vert_{{\cal R}} = (|a|^2 + |b|^2)^{1/2}\, ,
$$
hence by (\ref{a})  
$$
\lim_{\e \to 0}\sum_{i\in I_{\e}}\sum_{j\in J_{\e}} \e^3 | a_{\e,i}|^2 + 
{\e}^3 |b_{\e,i} - u_{\e,2,i,j}|^2 =0
$$
which implies that
\begin{eqnarray}\label{a1}
\lim_{\e \to 0}\sum_{i\in I_{\e}} \e^2 | a_{\e,i}|^2 =0
\end{eqnarray} 
and 
\begin{equation}\label{a12}
\sum_{i\in I_{\e}} \e^2 | b_{\e,i}|^2 \le c
\end{equation}
for each $\e>0$ small enough. 

Since $|x_j-x_i|\le \eta$, by the equivalence of the norms we have
\begin{eqnarray}\label{sommaai}
\nonumber
\sum_{i,j\in I_k}{1\over h^2}\int_{\e E^2_i} \Bigl| 
a_{\e,i} \wedge (x_j-x_i)\Bigr|^2 \dx 
&\le&
c \sum_{i,j\in I_k}{\e^2\over h^2}\eta^2 |a_{\e,i}|^2  \\
&=&
c\eta^2 \sum_{i\in I_k} \e^2 |a_{\e,i}|^2 \, .
\end{eqnarray}
Note that $\sum_{k\in J}\sum_{i\in I_k} =\sum_{i\in I_{\e}}$.

Now we insert
(\ref{sommaai}) into (\ref{beta}) and, summing up all the
corresponding estimates obtained for different indices $k\in J$, by
(\ref{trasl}) we get
\begin{eqnarray}\label{2term}
&&
\sum_{k\in J}\sum_{j\in I_k}\int_{\e E^2_j}\Bigl| 
 a_{\e,j} \wedge x +b_{\e,j}-
{1\over h^2}\sum_{i\in I_k} a_{\e,i} \wedge x +
b_{\e,i}\Bigr|^2\dx\\ \nonumber
&\le&
c\Bigl(\sum_{i\in I_{\e}}\int_{\e E^2_i} \Bigl|a_{\e,i} \wedge x +b_{\e,i}- 
u_2(x)\Bigr|^2\dx
+\eta^2 \sum_{i\in I_{\e}} \e^2 |a_{\e,i}|^2\\ \nonumber
&&\qquad
+\sup_{|\tau|\le \sqrt{2}\eta}\Vert u_2(\cdot)-u_2(\cdot +\tau)
\Vert^2_{{\rm L}^2(\Omega;\rr^3)}\Bigr) \,.
\end{eqnarray}

Finally, we sum up the estimates (\ref{gamma}) for $k\in J$ and insert 
(\ref{2term});
by (\ref{b1}) and (\ref{a1}) we get 
\begin{eqnarray}\label{a13}
&&
\lim_{\e\to0} \sum_{k\in J} 
\int_{T_{\eta}^k\cap\e E^2} \Bigl| u_2- A_{\e}^k\wedge x - B_{\e}^k \Bigr|^2
\dx \\ \nonumber
&\le&
c \sup_{|\tau|\le \sqrt{2}\eta}\Vert u_2(\cdot)-u_2(\cdot +\tau)
\Vert^2_{{\rm L}^2(\Omega;\rr^3)} \,.
\end{eqnarray}
On the other hand it is easy to see by (\ref{a1}) and (\ref{a12}) that there
exists $A^k\wedge x + B^k$ such that 
$$
\lim_{\e\to 0}\int_{T_{\eta}^k} \Bigl|A_{\e}^k\wedge x + B_{\e}^k -
A^k\wedge x - B^k \Bigr|^2 \dx =0
$$
for each $k\in J$,
hence by (\ref{a13}) we can conclude that $u_2\in {\cal U}$. 

Conversely, if $u_2\in {\cal U}$ then $\e E^2= \cup_{k\in J} 
T_{\eta}^k\cap \e E^2$ and we have rigid displacements $A^k\wedge x +B^k$ 
on each $T^k_{\eta}$.

We define
$$
a_{\e,i} \wedge x +b_{\e,i}= (A^k\wedge x+B^k)_{|\e E^2_i}
$$
for each $i\in I_k$.
Hence
$$
\sum_{i\in I_{\e}}\int_{\e E^2_i} \Bigl| 
a_{\e,i} \wedge x +b_{\e,i}- u_2(x)\Bigr|^2 \dx =
\sum_{k\in J} \int_{T_{\eta}^k\cap \e E^2}
\Bigl|A^k\wedge x +B^k - u_2(x)\Bigr|^2 \dx
$$
and by definition of ${\cal U}$
\begin{eqnarray}\label{2mem}
&&
\lim_{\e\to 0}\sum_{k\in J} \int_{T_{\eta}^k\cap \e E^2}
\Bigl|A^k\wedge x +B^k - u_2(x)\Bigr|^2 \dx \\ \nonumber
&=&
\sum_{k\in J}|E| \int_{T_{\eta}^k\cap \Om}
\Bigl|A^k\wedge x +B^k - u_2(x)\Bigr|^2 \dx 
\le o(\eta) \,.
\end{eqnarray}
By (\ref{2mem}), passing to the limit as $\eta\to0$, we get
$$
\lim_{\e\to0} \sum_{i\in I_{\e}}\int_{\e E^2_i} \Bigl| 
a_{\e,i} \wedge x +b_{\e,i}- u_2(x)\Bigr|^2 \dx = 0\, .
$$
\qed

\begin{remark}\label{rem}{\rm 
Note that, since $u_{\e}$ are rigid displacements, by (\ref{b2}) it is easy 
to see that $u_1$ is a rigid displacement. 

For simplicity, we will denote
$$
F(u_1,u_2;\Om)=\Gamma\hbox{-}\lim_{\e\to0} F^\gamma_\e(u_1,u_2;\Om)
$$
for $(u_1,u_2)\in {\cal R}\times {\cal U}$. 
We will continue to write $F^\gamma_\e(u,\Om)$ for $u\in
\ld^2_{\mu_\e}(\Om)$.
}
\end{remark}
\begin{theorem}\label{teorema}
For $\gamma=2$ the functionals $F^\gamma_\e$ $\Gamma$-converge as 
$\e \to 0$ to 
$$
F(u_1,u_2;\Om)= c_1 \int_{\Om} |(u_1)_{\alpha}-(u_2)_{\alpha}|^2 \dx + c_2 
\int_{\Om} |(u_1)_3 - (u_2)_3 |^2 \dx 
$$
on ${\cal R}\times {\cal U}$ 
with respect to the convergence introduced in Definition {\rm\ref{def}}, where 
$c_1= {3\over 8}\pi$, $c_2={\pi\over 4}$.
\end{theorem}
{\sc Proof.} By the invariance of the 
functionals with respect to translations of rigid displacements and by Remark 
\ref{rem} we can always assume without loss of 
generality that $u_{\e}=u_1$ on $\e E^1$.

Let us call
$$\alpha_{\e,i}\wedge x +\beta_{\e,i} = 
 u_1 - a_{\e,i} \wedge x - b_{\e,i} 
$$ 
hence
$$
F^\gamma_\e(u_{\e},\Om)=
\e^{\gamma-1}\sum_{i\in I_{\e}}\int_{\partial \e E_i^2}
\Bigl|(\alpha_{\e,i}\wedge x +
\beta_{\e,i}) \odot \nu  \Bigr|^2 d{\cal H}^2 \, .                 
$$

Fix $x_3\in (0,1)$, we can find the following equality 
\begin{eqnarray*}
\nonumber
&& 4 \e \int_{\partial \e D_i}\Bigl |(\alpha_{\e,i}\wedge x + \beta_{\e,i})
\odot \nu \Bigr|^2 d{\cal H}^1 -
16  \int_{\e D_i}\Bigl |\alpha_{\e,i}\wedge x + \beta_{\e,i}
\Bigr|^2 \dx_{\alpha}\\  \nonumber            
&=&
{\pi\over 2}\e^2 \Bigl( \Bigl|(\ave_{\e D_i}\alpha_{\e,i}\wedge x + 
\beta_{\e,i} \,\dx_{\alpha})_1 \Bigr|^2 + \Bigl|(\ave_{\e D_i}\alpha_{\e,i}
\wedge x + 
\beta_{\e,i} \,\dx_{\alpha})_2 \Bigr|^2 \Bigr) \\ 
&& +{\pi \over 64} \e^4 ((\alpha_{\e,i})_1^2+
(\alpha_{\e,i})_2^2+2 (\alpha_{\e,i})_3^2)\,. 
\end{eqnarray*}
Hence, if we integrate also in $x_3$, we get
\begin{eqnarray}\label{ug}
\nonumber
&& 4 \e \int_{\partial\e E^2_i}\Bigl |(\alpha_{\e,i}\wedge x + \beta_{\e,i})
\odot \nu \Bigr|^2 d{\cal H}^2 -
16  \int_{\e E^2_i}\Bigl |\alpha_{\e,i}\wedge x + \beta_{\e,i}
\Bigr|^2 \dx\\  \nonumber            
&=&
\int_0^1
{\pi\over 2}\e^2 \Bigl( \Bigl|(\ave_{\e D_i}\alpha_{\e,i}\wedge x + 
\beta_{\e,i} \,\dx_{\alpha})_1 \Bigr|^2 + \Bigl|(\ave_{\e D_i}\alpha_{\e,i}
\wedge x + 
\beta_{\e,i} \,\dx_{\alpha})_2 \Bigr|^2 \Bigr) \dx_3\\ 
&& +{\pi \over 64} \e^4 ((\alpha_{\e,i})_1^2+
(\alpha_{\e,i})_2^2+2 (\alpha_{\e,i})_3^2) \, .
\end{eqnarray}

But 
$$
\lim_{\e \to 0}\sum_{i\in I_{\e}} \int_{\e E^2_i} \Bigl|(\ave_{\e D_i}
\alpha_{\e,i} \wedge x + 
\beta_{\e,i}\,\dx_{\alpha})_h \Bigr|^2 \dx=
\lim_{\e \to 0} \sum_{i\in I_{\e}}\int_{\e E^2_i}\Bigl|(\alpha_{\e,i}\wedge x 
+ \beta_{\e,i})_h \Bigr|^2 \dx
$$
for each $h=1,2,3$, and 
$$
{\pi\over 2}\e^2 \Bigl|(\ave_{\e D_i}\alpha_{\e,i}\wedge x + 
\beta_{\e,i} \dx_{\alpha})_h \Bigr|^2   =
8 \int_{\e D_i} \Bigl|(\ave_{\e D_i}\alpha_{\e,i}\wedge x + 
\beta_{\e,i}\, \dx_{\alpha})_h \Bigr|^2 \dx_{\alpha}\, ;
$$
hence,
\begin{eqnarray}\label {c}
&&\nonumber
\lim_{\e \to 0} \sum_{i\in I_{\e}}\int_0^1{\pi\over 2}\e^2 \Bigl|(\ave_{\e D_i}
\alpha_{\e,i}\wedge x + 
\beta_{\e,i} \dx_{\alpha})_h\Bigr|^2 \dx_3\\
&=&
8 \lim_{\e \to 0} \sum_{i\in I_{\e}}\int_{\e E^2_i} \Bigl|(\alpha_{\e,i}\wedge
 x + 
\beta_{\e,i})_h \Bigr|^2 \dx\, .
\end{eqnarray}
If we pass to the limit in (\ref{ug}), by (\ref{c}) we obtain
\begin{eqnarray}\label{c1}
&&\nonumber
\liminf_{\e \to 0} \sum_{i\in I_{\e}}  \e \int_{\partial \e E^2_i}\Bigl 
|(\alpha_{\e,i}
\wedge x + \beta_{\e,i}) \odot \nu \Bigr|^2 d{\cal H}^2 \\
&\ge& \nonumber
6 \lim_{\e \to 0} \sum_{i\in I_{\e}} \int_{\e E^2_i} \Bigl|(\alpha_{\e,i}
\wedge x + 
\beta_{\e,i})_1 \Bigr|^2 + \Bigl|(\alpha_{\e,i}\wedge x + 
\beta_{\e,i})_2 \Bigr|^2 \dx\\
&&\nonumber
+ 4 \lim_{\e \to 0} \sum_{i\in I_{\e}} \int_{\e E^2_i}\Bigl|(\alpha_{\e,i}
\wedge x + 
\beta_{\e,i})_3 \Bigr|^2 \dx\\
&&
+ \lim_{\e \to 0} {\pi \over 64} \sum_{i\in I_{\e}}\e^4 ((\alpha_{\e,i})_1^2+
(\alpha_{\e,i})_2^2+2 (\alpha_{\e,i})_3^2)\, .
\end{eqnarray}

For every sequence $u_{\e}$ converging to $(u_1,u_2)$ in the sense of 
Definition \ref{def}, by (\ref{a1}) we have that
\begin{eqnarray}\label {d}
\lim_{\e \to 0} {\pi \over 64} \sum_{i\in I_{\e}}\e^4 ((\alpha_{\e,i})_1^2+
(\alpha_{\e,i})_2^2+2 (\alpha_{\e,i})_3^2) = 0 
\end{eqnarray}
so we insert (\ref{d}) into (\ref{c1}) to find that

\begin{eqnarray}\label{liminf}
\liminf_{\e \to 0} \e^2\int_{\Omega} \Bigl|{d Eu_{\e}\over d\mu_\e}\Bigr|^2 
d\mu_\e 
&\ge&
6 |E| \int_{\Om} |(u_1)_{\alpha}-(u_2)_{\alpha}|^2 \dx \\ \nonumber
&&
+ 4 |E| \int_{\Om} |(u_1)_3-(u_2)_3|^2 \dx\, .
\end{eqnarray}

By the arbitrarity of $u_{\e}$, choosing $\gamma=2$
\begin{equation}\label{gammaliminf}
\Gamma \hbox{-}\liminf_{\e\to 0} {F_{\e}}^2 (u_1,u_2;\Om)\ge 
F(u_1,u_2;\Om)\,.
\end{equation}
Now we consider  
$$
u_{\e}=  (c\wedge x + d) \,\chi_{ \e E^1} 
+ (a\wedge x + b)\, \chi_{ \e E^2}
$$ 
obviously it converges to $(c\wedge x + d, a\wedge x + b )$, 
and we call $\alpha \wedge x +\beta =(a-c)\wedge x +(b-d)$.

In this case  
$$
8 \int_{\e E^2_i}\Bigl|(\alpha\wedge x + \beta)_h \Bigr|^2 \dx
=
\int_0^1{\pi\over 2}\e^2 \Bigl|(\ave_{\e D_i}
\alpha\wedge x + 
\beta \dx_{\alpha})_h\Bigr|^2 \dx_3 +{\pi \over 128} {\e}^4 \alpha_3^2
$$
for $h=1,2$, hence by (\ref{ug})
\begin{eqnarray}\label{rigidi}
&&\nonumber
\limsup_{\e \to 0} \sum_{i\in I_{\e}}\e \int_{\partial \e E^2_i}\Bigl |(\alpha
\wedge x + \beta) \odot \nu \Bigr|^2 d{\cal H}^2 \\
&\le& \nonumber
6 \lim_{\e \to 0} \sum_{i\in I_{\e}} \int_{\e E^2_i} \Bigl|(\alpha\wedge x + 
\beta)_1 \Bigr|^2 + \Bigl|(\alpha\wedge x + 
\beta)_2 \Bigr|^2 \dx\\
&&\nonumber
+ 4 \lim_{\e \to 0} \sum_{i\in I_{\e}} \int_{\e E^2_i}\Bigl|(\alpha\wedge x + 
\beta)_3 \Bigr|^2 \dx
+ c \lim_{\e \to 0} \e^2 |\alpha|^2 \\
&=& 
6 |E| \int_{\Om} \Bigl|(\alpha\wedge x + 
\beta)_{\alpha} \Bigr|^2 \dx +4 |E| \int_{\Om} \Bigl|(\alpha\wedge x + 
\beta)_3 \Bigr|^2 \dx \,.
\end{eqnarray}
By (\ref{liminf}) and (\ref{rigidi}) we get
\begin{eqnarray}\label{limrigidi}
\lim_{\e\to0}\e^2\int_{\Omega} \Bigl|{d Eu_{\e}\over d\mu_\e}\Bigr|^2 
d\mu_\e &=&
6 |E| \int_{\Om} \Bigl|(\alpha\wedge x + 
\beta)_{\alpha} \Bigr|^2 \dx \\ \nonumber
&&
+4 |E| \int_{\Om} \Bigl|(\alpha\wedge x + 
\beta)_3 \Bigr|^2 \dx \,.
\end{eqnarray}
Now we fix $\eta>0$ and  consider $u_1\in {\cal R}$ and 
$v^{\eta}_2$ such that ${v^{\eta}_2}_{|T_{\eta}^k}= A^k\wedge x + B^k$ 
with $k\in J$. By (\ref{limrigidi}) we get
\begin{eqnarray}\label{limsupeta}
\nonumber
&&\limsup_{\e\to0} F^2_{\e}(u_1\,
\chi_{ \e E^1} + v^{\eta}_2\, \chi_{ \e E^2},\Omega)\\ \nonumber
&\le&
\sum_{k\in J} \limsup_{\e\to0} F^2_{\e}(u_1\,
\chi_{ \e E^1} +  (A^k\wedge x +B^k)\, \chi_{ \e E^2},
T_{\eta}^k\cap\Omega)\\ \nonumber
&=& \sum_{k\in J} 6 |E| \int_{T_{\eta}^k\cap\Omega} 
\Bigl|(u_1(x)- A^k\wedge x - B^k)_{\alpha}\Bigr|^2 \dx \\ \nonumber
&&
+ \sum_{k\in J} 4 |E| \int_{T_{\eta}^k\cap\Omega} 
\Bigl|(u_1(x)- A^k\wedge x - B^k)_3 \Bigr|^2 \dx \\ 
&=&
6 |E| \int_{\Omega} \Bigl|(u_1(x)- v^{\eta}_2(x))_{\alpha}\Bigr|^2 \dx 
+ 4 |E| \int_{\Omega} 
\Bigl|(u_1(x)- v^{\eta}_2(x))_3 \Bigr|^2 \dx
\end{eqnarray}
If $u_2\in {\cal U}$ then for all $\eta>0$ there exists $v^{\eta}_2$ as above 
such that
$\Vert u_2-v^{\eta}_2\Vert_{{\rm L}^2(\Omega;\rr^3)} \le o(\eta)$, since the
$\Gamma$-upper limit is ${\rm L}^2$-lower semicontinuous
if we denote
$$F''_2(u_1,u_2;\Omega)=
\Gamma \hbox{-}\limsup_{\e \to 0} {F_{\e}}^2 (u_1,u_2; \Omega)
$$
by (\ref{limsupeta}) we get
\begin{eqnarray*}
F''_2(u_1,u_2;\Omega)
&\le&
\liminf_{\eta\to0} F''_2(u_1,v^{\eta}_2;\Omega)\\
&\le&
\liminf_{\eta\to0} 6 |E| \int_{\Omega} 
\Bigl|(u_1(x)- v^{\eta}_2(x))_{\alpha}\Bigr|^2 \dx \\
&&
\qquad \quad 
+ 4 |E| \int_{\Omega} 
\Bigl|(u_1(x)- v^{\eta}_2(x))_3 \Bigr|^2 \dx \\
&=&
 6 |E| \int_{\Omega} \Bigl|(u_1(x)- u_2(x))_{\alpha}\Bigr|^2 \dx 
+ 4 |E| \int_{\Omega} 
\Bigl|(u_1(x)- u_2(x))_3 \Bigr|^2 \dx\,.
\end{eqnarray*}

It follows that given $(u_1,u_2)\in {\cal R}\times {\cal U}$
$$
\Gamma \hbox{-}\limsup_{\e \to 0} {F_{\e}}^2 (u_1,u_2;\Om)\le F(u_1,u_2;\Om)
$$
so that by (\ref{gammaliminf})
$$
\Gamma \hbox{-}\lim_{\e \to 0} {F_{\e}}^2 (u_1,u_2;\Om)= F(u_1,u_2;\Om)
$$
as desired.
\qed

If $u_{\e}$ converges to $(u_1,u_2)$ in the sense of Definition \ref{def} 
then $u_{\e}$
converges weakly in ${\rm L}^2(\Om;\rr^3)$ to $(1-c) u_1 + c u_2$ where 
$c=|E|$. If we define the energy
$$
F(u,\Om):= \inf_{\stackrel{\scriptstyle u=(1-c) u_1 + c u_2}
{(u_1,u_2)\in {\cal R}\times {\cal U}}}
F(u_1,u_2;\Om)
$$
by Theorem \ref{teorema} 
$$
F(u,\Om) = \inf_{r\in {\cal R}} \Bigl(\widetilde{c_1} \int_{\Om} |r_{\alpha} - 
u_{\alpha}|^2 \dx
+ \widetilde {c_2} \int_{\Om} |r_3 - u_3|^2 \dx\Bigr) 
$$
where $\widetilde{c_1}=c_1/c^2$ and $\widetilde{c_2}=c_2/c^2$, which
explains the non local nature of our limit.

\begin{remark}\label{}{\rm 
Let us consider, up to normalization, the same measure of 
Example \ref{24.9bis} 
$$
\widetilde{\mu_\e}(B)= \Bigl(|B| + \e{\cal H}^2 (B\cap \partial \e E^2)
\Bigr) 
$$
and the functionals
$$
\widetilde{F^2_{\e}}(u,\Om)=\e^2\int_{\Omega} \Bigl|{d Eu_{\e}\over 
d\widetilde{\mu_\e}}\Bigr|^2 d\widetilde{\mu}_\e \, .
$$
In this case by Theorem \ref{teorema} we can deduce that the 
$\Gamma\hbox{-}\limsup_{\e\to 0} \widetilde{F^2_{\e}}(u_1,u_2;\Om)$ is 
finite for $(u_1,u_2) \in {\cal R}\times {\cal U}$.

In fact, since
$\ld^2_{\mu_\e}(\Om;\rr^3)\subset \ld^2_{\widetilde{\mu_\e}}
(\Om;\rr^3) $,
given $(u_1,u_2)\in {\cal R}\times {\cal U}$  we have
$$
\Gamma\hbox{-}\limsup_{\e\to 0} \widetilde{F^2_{\e}}(u_1,u_2;\Om)\le
\Gamma\hbox{-}\limsup_{\e\to 0} F^2_{\e}(u_1,u_2;\Om)\,.
$$

}\end{remark}

\smallskip
\noindent{\bf Acknowledgements}\ \ We wish to express our thanks to 
Prof. Andrea Braides for suggesting the problem and 
for many stimulating conversations. We also thank Prof. Luigi Ambrosio 
for helpful comments.

\vspace{2cm}
\begin{tabular}{ll}
\qquad  Nadia Ansini\qquad\qquad & Fran\c{c}ois Bille Ebobisse \\
\qquad  SISSA/ISAS\qquad\qquad   & Dipartimento di Matematica \\
\qquad  Via Beirut 4, 34014 Trieste, Italy\qquad\qquad & 
Universit\`a di Pisa \\
 & Via Buonarroti 2, 56127 Pisa, Italy  
\end{tabular}

\end{document}